\documentclass[12pt]{article}
\input{psfig}
 \usepackage{amssymb}
 
 \newcommand{\C}{{\mathbb C}}
 \newcommand{\Z}{{\mathbb Z}}
 \newtheorem{theorem}{Theorem}
 \newtheorem{lemma}{Lemma}

 \def\Box
  {\thinspace\vbox{\hrule height .5pt \hbox{\vrule
   width .5pt \vbox to 7pt{\hbox to 3.5pt{}} \vrule width .5pt}
   \hrule height 0pt depth .5pt}}

 \newenvironment{proof}{{\it Proof:\/}}{$\Box$\vskip 0.08in}
 \newcommand{\mod}{{\mbox{ mod }}}
 
 \newcommand{\wt}{\widehat}
 \newcommand{\diag}[2]{\parbox{#2}{\psfig{figure=skein#1.eps,height=#2}}}

\begin{document}
\title{Skein modules at the $4$th roots of unity}
\author{Adam S. Sikora\\
{\em \small Department of Mathematics, University of Maryland,}\\
{\em \small College Park, MD 20742}\\
{\em \small asikora@math.umd.edu}}
\date{}

\maketitle

\begin{abstract}
The Kauffman bracket skein modules, ${\cal S}(M, A),$ have been
calculated for $A=\pm 1$ for all $3$-manifolds $M$ by relating them to the
$SL_2(\C)$-character varieties.
We extend this description to the case when $A$ is a
$4$th root of $1$ and $M$ is either a ${\rm surface}\times [0,1]$
or a rational homology sphere (or its submanifold).
\end{abstract}

\section{Introduction}
Let $R$ be a fixed commutative ring with $1$ and with a specified
invertible element $A\in R.$ Given an oriented $3$-manifold $M$
(possibly with boundary), the set of unoriented, framed links in
$M,$ considered up to isotopy, is denoted by ${\cal L}(M).$ 
The ``empty link,'' $\emptyset,$ is also an element of 
${\cal L}(M).$ The
Kauffman bracket skein module\footnote{We will also call it
{\em skein module} for short.} of $M$ is the
quotient ${\cal S}(M, R, A)=R{\cal L}(M)/{\cal S},$ where $R{\cal
L}(M)$ is the free $R$-module spanned by ${\cal L}(M)$ and ${\cal
S}$ is its submodule generated by the Kauffman bracket skein
relations:
\begin{equation}
\diag{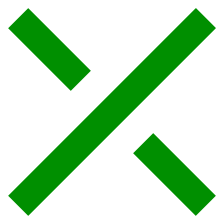}{.6cm}-A\diag{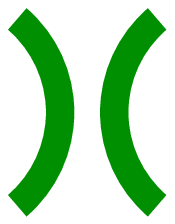}{.6cm}-A^{-1}\diag{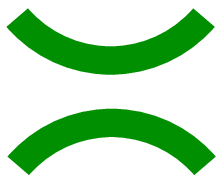}{.6cm},\label{e1}
\end{equation}
\begin{equation}
L\cup \diag{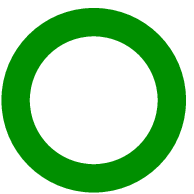}{.5cm} + (A^{2}+A^{-2})L.\label{e2}
\end{equation}
Above, $L$ is any link unlinked with the trivial knot, $\diag{0}{.5cm}.$
For more about skein modules see \cite{P, P-S-2, BFK} and the
bibliography therein.

Despite their simple definition, skein modules are much more
difficult to handle than it may appear at first sight. For
example, no algorithm is known for calculating the basis of ${\cal
S}(M,K,A)$ for a given $3$-manifold $M,$ field $K$ and $A\in K^*.$
On the other hand, skein modules are well understood for any ring
of coefficients $R$ and $A=\pm 1,$ since in this case a
particularly simple skein relation is satisfied:
$$\diag{3}{.6cm}=\diag{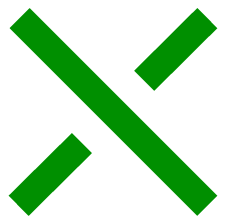}{.6cm}.$$ Thus, 
for any two links $L_1,L_2\subset M$ their disjoint
union, $L_1\cup L_2$ is a uniquely defined
element of ${\cal S}(M,R,A).$ Hence ${\cal S}(M,R,A)$ for $A=\pm 1$ becomes a 
commutative $R$-algebra, with the product of $L_1,L_2\subset M$ given 
by $L_1\cup L_2\in {\cal S}(M,R,A).$ This algebra is isomorphic to the
{\it skein algebra of $\pi_1(M)$} introduced in
\cite{P-S-1}:\vspace{.15in}

\noindent{\bf Definition} {\em Given a commutative ring $R$ with
$1$ and a group $G$ with the identity $e\in G,$ the skein algebra
of $G$ is the quotient, $${\cal S}(G,R)=R[x_g\ g\in G]/I,$$ where
$R[x_g\ g\in G]$ is the ring of polynomials in variables indexed by elements
of $G$ and $I$ is the ideal generated by $x_e-2$ and
$x_gx_h-x_{gh}-x_{gh^{-1}}$ for any $g,h\in G.$}

Assume that $N$ is a rational homology sphere with no
$2$-torsion in $H_1(N,\Z),$ and $M$ is a submanifold of $N.$ Let
$A$ be a primitive $4$th root of $1$ in $R.$ Our first result
states that ${\cal S}(M,R,A)$ is a commutative and associative
$R$-algebra with the product of any two links $L_1,L_2 \in {\cal
S}(M,R,A)$ given by $(-1)^{d(L_1,L_2)}L_1\cup L_2\in {\cal
S}(M,R,A),$ where $d(L_1,L_2)\in \Z/4\Z.$ The definition of
$d(L_1,L_2)$ is somewhat complicated, and therefore we
precede it by a few remarks and a lemma.

If $r$ is the smallest positive integer annihilating the torsion
part of $H_1(N,\Z),$ $rTH_1(N,\Z)=0,$ then the linking number
$lk(L_1,L_2),$ defined for any two disjoint oriented links
$L_1,L_2\subset M,$ takes values in $\frac{1}{r}\Z.$ Since the
fractional part of $lk(L_1,L_2)$ depends on $[L_1],[L_2]\in
H_1(N,\Z)$ only, there is a well defined linking form $lf:
H_1(N,\Z)\times H_1(N,\Z)\to \frac{1}{r}\Z/\Z,$
$$lf([L_1],[L_2])=lk(L_1,L_2) \quad{\rm mod\ } \Z.$$

\begin{lemma}\label{l1} (Proof in Section \ref{sproofs})
If $H=H_1(N,\Z)$ has no $2$-torsion part then the linking form,
$lf: H\times H\to \frac{1}{r}\Z/\Z,$ can be lifted to a
symmetric bilinear form $\wt{lf}: H\times H\to
\frac{1}{r}\Z/4\Z.$
\end{lemma}

For any two disjoint oriented links $L_1,L_2$ in $M$ 
we define 
$$d(L_1,L_2)=\wt{lf}([L_1],[L_2])-lk(L_1,L_2).$$

Since
$$\wt{lf}([L_1],[L_2])-lf([L_1],[L_2])\in \Z/4\Z {\rm \ and\ }
lf([L_1],[L_2])-lk(L_1,L_2)\in \Z,$$ $d(L_1,L_2)$ takes values in $\Z/4\Z.$
Observe that $d(L_1,L_2)$ mod $2$ does not depend on orientations
of $L_1$ and $L_2.$

\begin{theorem} \label{th1} (Proof in Section \ref{sproofs})
If $M$ is a submanifold of a rational homology sphere $N$ without
$2$-torsion in $H_1(N,\Z)$ then
${\cal S}(M,R,A)$ is a commutative, associative $R$-algebra with
the multiplication defined by
\begin{equation}
L_1\cdot L_2=L_1\cup L_2\cdot (-1)^{d(L_1,L_2)},\label{e-product}
\end{equation}
for any two disjoint framed, unoriented links in $L_1,L_2\subset
M.$
\end{theorem}

\begin{theorem}\label{th2}(Proof in Section \ref{sproofs})
Under the assumptions of Theorem \ref{th1}, there is an isomorphism of
$R$-algebras $$\Psi: {\cal S}(\pi_1(M),R)\to {\cal S}(M,R, A),$$
such that $\Psi(x_g)=K\cdot
A^{D(K)},$ where $K$ is an arbitrary framed unoriented knot which
represents the conjugacy class of $g^{\pm 1}\in
\pi_1(M)$ and $D(K)\in \Z/4\Z.$ \footnote{$D(K)$ will be defined in
Section \ref{sproofs}.}
\end{theorem}

Since Lemma \ref{l1} does not hold for groups $H$ with a 
$2$-torsion, we do not know how to extend the statements of Theorems
\ref{th1} and \ref{th2} to submanifolds of all rational homology spheres.
On the other hand Theorems
\ref{th1} and \ref{th2} cannot be easily extended to all $3$-manifolds.
If $M$ is not a submanifold of a rational homology sphere then
$M$ contains a non-separating surface. It is known, that if such
surface is a torus then ${\cal S}(M,R,\pm 1)$ is not isomorphic
to ${\cal S}(M,R,A)$ for $A^4=1,$ $A\ne \pm 1.$
On the other hand, Theorems \ref{th1} and \ref{th2} imply that
$${\cal S}(M,R,A)\simeq {\cal S}(\pi_1(M),R)={\cal S}(M,R,\pm 1),$$
for any $4$th roof of $1$, $A.$

\section{Skein algebras of surfaces}
The third result  of this paper concerns skein algebras of
surfaces. If $F$ is an oriented surface then the skein module of
$F\times [0,1]$ has a product structure for any $R$ and $A:$ the
product of two links $L_1,L_2\subset F\times [0,1]$ is a union of
them, $L_1\cup L_2,$ such that $L_1$ lies in the upper half of the
cylinder and $L_2$ lies in the lower part. For that reason, ${\cal
S}(F\times [0,1],R,A)$ is called the skein algebra of $F,$ and it
is denoted by ${\cal S}(F,R,A).$ \footnote{The reason for that
notation is that the algebra structure on ${\cal S}(M,R,A)$ for
$M=F\times [0,1]$ is not uniquely determined by $M:$ if $F_1$ is a
punctured torus and $F_2$ is a twice punctured disk then
$F_1\times [0,1]$ and $F_2\times [0,1]$ are homeomorphic but 
${\cal S}(F_1\times [0,1],R,A)$ and ${\cal S}(F_2\times [0,1],R,A)$ 
are not isomorphic as algebras.} Note, that if $A\ne \pm 1$ then 
this product will be
usually non-commutative and therefore different from the products
considered before. The skein algebras of surfaces have been
described for a few surfaces in \cite{B-P, F-G, P-S-2, S}. An
interesting connection between the skein algebra of a torus and
non-commutative geometry has been discovered in \cite{F-G}.
Unfortunately, a complete description of skein algebras of
surfaces is known for $A=\pm 1$ only: ${\cal S}(F,R,\pm 1)\simeq 
{\cal S}(\pi_1(F),R).$ We extend
this description to all $4$th roots of $1,$ $A,$ by showing that
${\cal S}(F,R,A)$ is isomorphic to the skein algebra of $\pi_1(F)$
deformed along the symplectic intersection form $\omega:
H_1(F,\Z)\times H_1(F,\Z)\to \Z.$\vspace{.15in}

\noindent{\bf Definition} {\em Let $R<x_g,\ g\in G>$ be the free
associative, non-commutative $R$-algebra with $1$ whose (free)
generators, $x_g,$ are in a bijection with the elements of $G.$ If
$\omega$ is a bilinear form on $H=G/[G,G]$ and $A$ is an invertible
element in $R,$ then {\it the skein
algebra of $G$ deformed with respect to $\omega$} is $${\cal
S}_{\omega}(G,R,A)=R<x_g,\ g\in G>/I,$$ where $I$ is the ideal
generated by $x_e-2$ ($e$ is the identity in $G$) and expressions
of the form
\begin{equation}
x_gx_h- A^{\omega(g,h)}x_{gh}-
A^{-\omega(g,h)}x_{gh^{-1}},\label{e3}
\end{equation}
\begin{equation}
x_{hgh^{-1}}=x_g,\label{e3.5}
\end{equation}
for all $g,h\in G.$
$\omega(g,h)$ means $\omega([g],[h]),$ where $[g],[h]$ are classes of
$g,h$ in $H=G/[G,G].$}

It is easy to prove (see \cite[Fact 2.6]{P-S-1})
that $x_{hgh^{-1}}=x_g$ in ${\cal S}(G,R)$ for any
$g,h\in G.$ Therefore, we have an isomorphism of algebras
\begin{equation}
{\cal S}_{0}(G,R,A)={\cal S}(G,R).\label{e-omega0}
\end{equation}

The definition of ${\cal S}_{\omega}(G,R,A)$ makes sense if either 
$\omega$ is an
integer-valued bilinear form on $H$ and $A$ is an arbitrary invertible element
of $R,$ or if $A$ is an $n$-th root of $1$ and $\omega$ assumes values
in $\Z/n\Z.$ In this paper, we assume that $A$ is a $4$th root of unity, and
$\omega$ is a skew-symmetric form. We will prove in Lemma \ref{prop_of_S(G)}
that under these assumptions, $x_gx_h=\pm x_hx_g$ for any $g,h\in G.$
Hence ${\cal S}_{\omega}(G,R,A)$ is similar to super-symmetric algebras.

\begin{theorem}\label{th3}
For any oriented surface $F$ there is an isomorphism of
$R$-algebras $$\Psi: {\cal S}_{\omega}(\pi_1(F),R,A)\to {\cal
S}(F,R, A),$$ such that $\Psi(x_g)=K\cdot A^{D(K)},$ where $K$ is
an arbitrary framed unoriented knot which represents the conjugacy
class of $g^{\pm 1}\in \pi_1(F).$\footnote{The symbol $D(K)\in
\Z/4\Z$ will be defined in Section \ref{sproofs}.} Here, $\omega$ denotes the
standard symplectic $2$-form on $H_1(F,\Z).$
\end{theorem}

Note that the above theorem is not a special case of Theorem \ref{th2} --
the algebras considered in Theorem \ref{th3} are non-commutative.
It is an open and interesting problem to determine if the above result can be
extended to other values of $A.$ This certainly can be done for tori --
it is easy to prove using the ``product-to-sum'' formula of \cite{F-G} 
(see also \cite{S}) that the skein algebra of a torus is isomorphic to ${\cal
S}_{\omega}(\Z\times \Z,R,A)$ for any $A.$

\section{The proofs}
\label{sproofs}

We assume as before that $N$ is a rational homology sphere,
and $r$ a minimal positive integer annihilating the torsion part of
$H_1(N,\Z),$  $rTH_1(N,\Z)=0.$ Moreover, we assume that $H_1(N,\Z)$
has no $2$-torsion and, hence, $r$ is odd. Let $M$ be a $3$-dimensional
submanifold of $N,$ and $R$ a commutative ring with a fixed $4$th
primitive root of $1,$ denoted by $A.$

\noindent{\bf Proof of Lemma \ref{l1}:} $H=H_1(N,\Z)$ decomposes into a sum
of finite cyclic groups of orders $n_1,n_2,...,n_k,$ whose
generators we denote by $g_1,...g_k.$ Consider a $k\times k$
matrix $(a_{ij})\in M_k(\frac{1}{r}\Z/4\Z)$ defined as follows:
Given $i,j,$ write $lf(g_i,g_j)$ as $\frac{k}{l},$ where $k,l$ are
relatively prime and $k\in\{0,1,...,l-1\}.$ Since $r$ is odd, $l =
\pm 1$ mod $4$ and, therefore, there is a unique
$\delta\in\{0,1,2,3\}$ such that $k+\delta l=0$ mod $4.$ Define
$a_{ij}$ to be $\frac{k}{l}+\delta\in \frac{1}{r}\Z/4\Z.$ The
matrix $(a_{ij})\in M_k(\frac{1}{r}\Z/4\Z)$ determines a bilinear
form $H\times H\to \frac{1}{r}\Z/4\Z$ if and only if 
$n_ia_{ij}=0$ mod $4$ and $n_ja_{ij}=0$ mod $4$ for any $i,j.$
We will show that $(a_{ij})$ indeed satisfies
this condition. Choose $i,j\in\{1,...,k\}$ and
denote $lf(g_i,g_j)$ by $k/l$ as before. Since
$n_ilf(g_i,g_j)=lf(n_ig_i,g_j)=0$ mod $\Z$ and, similarly,
$n_jlf(g_i,g_j)=0$ mod $\Z,$ $n_i$ and $n_j$ are divisible by $l.$
Hence $n_ia_{ij}=\frac{n_i}{l}(k+\delta l)=0$ mod $4.$ Similarly,
$n_ja_{ij}=0$ mod $4.$ Therefore $(a_{ij})$ determines a bilinear
form $\wt{lf}:H\times H\to \frac{1}{r}\Z/4\Z.$ For $x=\sum
\alpha_ig_i,\ y=\sum \beta_jg_j\in H,$
$$\wt{lf}(x,y)=\sum\alpha_i\beta_ja_{ij}.$$ Since $a_{ij}=a_{ji},$
$\wt{lk}$ is symmetric. \Box

So far we have defined the functions $lk, lf, \wt{lf}$ for pairs
of links in rational homology spheres. Now, we are going to define
analogous functions, denoted by the same symbols, for pairs of
links in cylinders over surfaces, $F\times [0,1].$ These functions
will have exactly the same properties as the functions defined
before, and because of that, it will be possible to prove Theorems
\ref{th2} and \ref{th3} together (as Theorem \ref{main}).

Let $F$ be an oriented surface (with or without boundary) and let
let $L_1,L_2$ be oriented links in $F\times [0,1]$ whose
projections on $F$ are in general position. Define $lk(L_1,L_2)$
to be $$\frac{1}{2}\sum signs\ of\ crossings\ between\ L_1\ and\
L_2\ \in \frac{1}{2}\Z.$$ Let $H=H_1(F,\Z)$ 
and let $lf:H\times H\to \frac{1}{2}\Z/\Z$ be the bilinear form
$$lf(c_1,c_2)=lk(L_1,L_2)\mod \Z,$$  where $L_1,L_2$ are any
oriented links in $F\times [0,1]$ representing the cycles
$c_1,c_2\in H.$ There is a symplectic form on $H:$
$$\omega:\ H\times H\to \Z,$$ $\omega(c_1,c_2)=\sum\ types\ of\
crossings\ between\ c_1\ and\ c_2.$\\
\centerline{{\psfig{figure=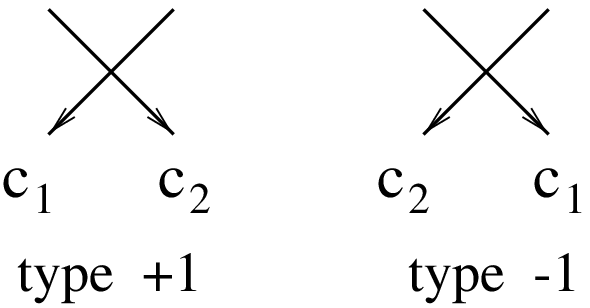,height=2cm}}} Observe that
$lf(c_1,c_2)=\frac{1}{2}\omega(c_1,c_2)\mod \Z$ and therefore $lf$
can be lifted to $\wt{lf}:\ H\times H\to \frac{1}{2}\Z/4\Z,$
$\wt{lf}=\frac{1}{2}\omega\mod 4\Z.$ Let
$d(L_1,L_2)=\wt{lf}([L_1],[L_2])-lk(L_1,L_2)\in \Z/ 4\Z.$ Note
that $d(L_1,L_2)$ is equal (mod $4$) to the number of
intersections between $L_1$ and $L_2$ in which $L_2$ lies on top
of $L_1.$  Exactly as before, $d(L_1,L_2)\mod 2\Z$ does not depend
on the orientations of $L_1$ and $L_2.$ This time however, $d$ is
not always symmetric. 

Let us summarize the notation introduced so far and fix it for the the
remainder of the paper. We will always assume that one of the two
following conditions holds:

\begin{enumerate}
\item[(I)] $M$ is a submanifold of a rational homology sphere; or
\item[(II)] $M=F\times [0,1],$ where $F$ is a surface. 
\end{enumerate}

In either case, we will assume that $lk,lf,\wt{lf}$ and $d$ are defined
appropriately as before, separately for (I) and (II).
$H$ denotes $H_1(N,\Z)$ in case (I), and $H_1(F,\Z)$ in case (II).
We define $\omega:\ H\times H\to \Z$ to be $0$ if (I) holds or to be
the symplectic form (defined above) if (II) holds.

Several basic properties of $d$ are stated below, common for
both cases, (I) and (II). The first of them is the definition of
$d.$ The proof of the remaining ones is left to the reader.

\begin{lemma}\label{prop_of_d}
\begin{enumerate}
\item[(1)] $d(L_1,L_2)=\wt{lf}([L_1],[L_2])-lk(L_1,L_2)\in \Z/4\Z$
\item[(2)] $d(L_1,L_2)\in \Z/2\Z$ does not depend on the orientations of
$L_1,L_2.$
\item[(3)] $d(L_1,L_2\cup L_3)=d(L_1,L_2)+d(L_1,L_3)$ and
$d(L_1\cup L_2,L_3)=d(L_1\cup L_3)+d(L_2,L_3)$
\item[(4)] If $L_1,L_2$ are two disjoint links in $M,$  and $L_1'$ is
obtained from $L_1$ by a change of crossing with $L_2$ from a negative
one to a positive one,
$\diag{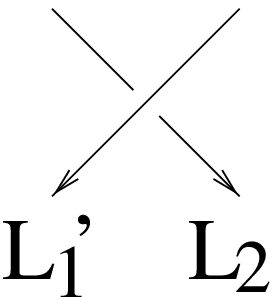}{1cm},\ \diag{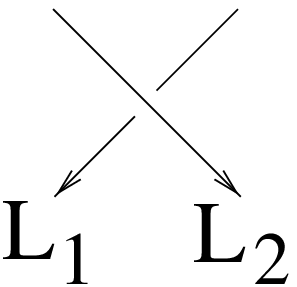}{1cm},$
then $d(L_1',L_2)= d(L_1,L_2)+1.$
\item[(5)] For any three skein related unoriented links,
\diag{3}{.6cm},\diag{4}{.6cm},\diag{5}{.4cm}, 
and for any unoriented link $L$ disjoint with them,
$d(L,\diag{3}{.6cm})=d(L,\diag{4}{.6cm})=d(L,\diag{5}{.6cm})\in
\Z/2\Z.$ \footnote{By property (2) above, this statement
makes sense.}
\item[(6)] $d(\diag{0}{.5cm},L)=0$ for any link $L$ unlinked with
the trivial knot, $\diag{0}{.5cm}.$
\item[(7)] For any two oriented knots $K_1,K_2,$
$d(K_1,K_2)-d(K_2,K_1)= \omega([K_1],[K_2])\in
\Z/4\Z.$
\end{enumerate}
\end{lemma}
\Box

The proof of Theorem \ref{th1} will be based on the fact that the
skein modules for $A$ being a primitive $4$th root of $1$ satisfy a
particularly simple skein relation:
\begin{equation}
\diag{3}{.6cm}=-\diag{3a}{.6cm}\label{e4}
\end{equation}
Indeed, the $90^o$ rotation of the diagrams in the skein relation
$$\diag{3}{.6cm}=A\diag{4}{.6cm}+A^{-1}\diag{5}{.6cm},$$
together with the substitution $A^{-1}=-A,$ yields
$$\diag{3}{.6cm}=A\diag{5}{.6cm}+A^{-1}\diag{4}{.6cm}=
-A^{-1}\diag{5}{.6cm}-A\diag{4}{.6cm}=-\diag{3a}{.6cm}.$$

By (\ref{e4}), the product of links $L_1, L_2\subset 
F\times [0,1]$ in ${\cal S}(F,R,A)$ is equal to
$L_1\cup L_2$ or $-L_1\cup L_2$ depending if the number of
crossings in which $L_2$ goes over $L_1$ is even or odd. 
We have noticed before, that this number is equal to
$d(L_1,L_2)$ mod $4$.
Therefore the product in skein algebras of surfaces is
given by
$$L_1\cdot L_2= (-1)^{d(L_1,L_2)}L_1\cup L_2.$$
Theorem \ref{th1} states that this formula defines a meaningful product
on ${\cal S}(M,R,A)$ in case (I) as well.

\begin{lemma}\label{l2}
The expression
\begin{equation} (-1)^{d(L_1,L_2)} L_1\cup L_2\in
{\cal S}(M,R,A)\label{e-expression}
\end{equation}
 depends on the isotopy types of $L_1$ and $L_2$
only, ie. it does not depend on the position of $L_1$ with respect
to $L_2.$
\end{lemma}

\begin{proof}
We need to prove that if $L_i$ and $L_i'$ are isotopic for $i=1,2$
then $L_1\cdot L_2=L_1'\cdot L_2'.$ We can assume that $L_1\cup
L_2$ and $L_1'\cup L_2'$ differ by one crossing change
only (and repeat the argument as many times as necessary). By
(\ref{e4}) and by Lemma \ref{prop_of_d}(4) $L_1\cup L_2= -L_1'\cup
L_2'$ and $d(L_1,L_2)=d(L_1',L_2')+1 \mod 2.$ Therefore
$$(-1)^{d(L_1,L_2)} L_1\cup L_2=(-1)^{d(L_1',L_2')}
L_1'\cup L_2'.$$
\end{proof}

{\bf Proof of Theorem \ref{th1}:}
By the lemma above, there is a bilinear form
$R{\cal L}(M)\times R{\cal L}(M)\stackrel{\cdot}{\to} {\cal
S}(M,R,A),$ $$L_1\cdot L_2= (-1)^{d(L_1,L_2)} L_1\cup L_2\in {\cal
S}(M,R,A).$$ We are going to show that this form descends to a
product on ${\cal S}(M,R,A),$ that is, that it vanishes on ${\cal
S}\times R{\cal L}(M)$ and $R{\cal L}(M)\times {\cal S},$ where
${\cal S}$ is the submodule of $R{\cal L}(M)$ generated by the
skein relations, (\ref{e1}) and (\ref{e2}). Since the bilinear
form $\cdot$ is symmetric, it is sufficient to show that $R{\cal
L}(M)\cdot {\cal S}=0,$ ie. that for any link $L$
$$L\cdot\diag{3}{.6cm}=AL\cdot\diag{4}{.6cm}+A^{-1}L\cdot\diag{5}{.5cm}
\in {\cal S}(M,R,A)$$ and
$$L\cdot(L'\cup \diag{0}{.5cm})=L\cdot 2L'.$$
The first identity holds since by Lemma \ref{prop_of_d}(5) it can
be rewritten as
$$L\cup\diag{3}{.6cm}=AL\cup\diag{4}{.6cm}+A^{-1}L\cup\diag{5}{.5cm}
\in {\cal S}(M,R,A).$$ The second identity is equivalent to
\begin{equation}
(-1)^{d(L,L'\cup \diag{0}{.4cm})}L\cup L'\cup \diag{0}{.5cm}
=2(-1)^{d(L,L')}L\cup L'. \label{e-skein2}
\end{equation}

Since $\diag{0}{.5cm}$ is the trivial knot, we can unlink $L$ and
$\diag{0}{.5cm}$ by a sequence of crossing changes between
$\diag{0}{.5cm}$ and $L.$ Such crossing changes will not change
the right nor the left side of the above identity (by Lemma
\ref{l2}).
Therefore it is enough to prove (\ref{e-skein2}) under the assumption
that $\diag{0}{.5cm}$ and $L$ are unlinked. In this situation we have
$$L\cup L'\cup
\diag{0}{.5cm}=2L\cup L'$$ and, by Lemma \ref{prop_of_d}(3) and
(6), $d(L,L'\cup \diag{0}{.5cm})=d(L,L').$ Hence identity
(\ref{e-skein2}) holds and, therefore, (\ref{e-product}) defines a
product on ${\cal S}(M,R,A).$ Since $d(L_1,L_2)=d(L_2,L_1),$
the product is commutative.
By Lemma \ref{prop_of_d}(3), we have
$$(L_1\cdot L_2)\cdot L_3=(-1)^{d(L_1,L_2)+d(L_1\cup L_2,L_3)}
L_1\cup L_2\cup L_3=$$
$$(-1)^{d(L_2,L_3)+d(L_1, L_2\cup L_3)} L_1\cup L_2\cup L_3=
L_1\cdot(L_2\cdot L_3),$$ and therefore $\cdot$ is associative.
\Box

The remainder of this paper is devoted to proving Theorems
\ref{th2} and \ref{th3}. Their statements use a function $D,$
which will be defined now. Let $M$
satisfy either condition (I) or (II). Let $L$ be a framed oriented
link in $M$ composed of $n$ connected components, $K_1,K_2,...,K_n.$ We
say that $\partial L=\partial_+ L\cup
\partial_- L$ is a splitting of $\partial L$ if $\partial_+ L$ is
a union of $n$ connected components of $\partial L,$ composed of
exactly one component from each $\partial K_i,$ and $\partial_- L$
is the union of the remaining $n$ connected components of
$\partial L.$ \footnote{Each connected component of $L$ is an
annulus embedded in $M,$ and hence it has two boundary
components.}

\begin{lemma}\label{l3}
$d(\partial_+ L,\partial_- L)$ does not depend on the choice of a splitting
$\partial L=\partial_+ L\cup \partial_-L.$
\end{lemma}

\begin{proof}
Assume first that $K$ is a framed oriented knot and $\partial K=\partial_+ K
\cup \partial_- K$ is one of its two possible splittings.
We need to prove that
\begin{equation}\label{e-D-for-knot}
d(\partial_+ K,\partial_-K)=d(\partial_- K,
\partial_+K).
\end{equation}
This holds in case (I) since then $d(\cdot,\cdot)$ is symmetric.
In case (II) we have $$d(\partial_+ K,\partial_-K)=
\wt{lf}(\partial_+ K,\partial_-K)-lk(\partial_+ K,\partial_-K)=
-lk(\partial_+ K,\partial_-K),$$ because $[\partial_+
K]=[\partial_-K]$ in $H$ and $\wt{lf}(x,x)=0$ for any $x\in H.$
Now, equality (\ref{e-D-for-knot}) follows from the fact that $lk$
is symmetric.

Assume now that $L$ has $n$ connected components, $K_1,...,K_n.$
In order to establish Lemma \ref{l3} it is enough to prove that
$d(\partial_+ L,\partial_-L)$ will not change if we change the
splitting $\partial L=\partial_- L\cup\partial_+L$ by one
connected component. Let $L'=L\setminus K_s,$ for some
$s\in\{1,2,...,n\},$ and let the new splitting be $L=(\partial_+
L'\cup \partial_-K_s)\cup (\partial_- L'\cup \partial_+ K_s).$ We
need to show that
\begin{equation}
d(\partial_+ L'\cup \partial_-K_s,\partial_- L'\cup \partial_+
K_s)=d(\partial_+ L'\cup \partial_+ K_s,\partial_- L'\cup
\partial_- K_s). \label{e-splitting-change}
\end{equation}
By Lemma \ref{prop_of_d}(3), the left side equals $$d(\partial_+
L',\partial_- L')+ d(\partial_+ L',\partial_+ K_s)+ d(\partial_-
K_s,\partial_- L')+d(\partial_- K_s,\partial_+ K_s),$$ and the
right side is equal to $$d(\partial_+ L',\partial_- L')+
d(\partial_+ L',\partial_- K_s)+ d(\partial_+ K_s,\partial_-
L')+d(\partial_+ K_s,\partial_- K_s).$$ Note that $\partial_+ K_s$
is isotopic to $\partial_- K_s$ in $M\setminus L',$ and hence
$d(\partial_+ L', \partial_+ K_s)=d(\partial_+ L', \partial_-
K_s).$ Similarly, $d(\partial_-K_s,\partial_-L')=d(\partial_+
K_s,\partial_- L').$ Therefore, by (\ref{e-D-for-knot}), the right
and the left side of (\ref{e-splitting-change}) are equal.
\end{proof}

Since $d(\partial_+ L,\partial_- L)$ does not depend on the
splitting, we denote it shortly by $D(L)\in \Z/4\Z.$

\begin{lemma}\label{d_no_orient}
$D(K)$ does not depend on the orientation of $K.$
\end{lemma}

\begin{proof}
If $K^-$ is the knot $K$ with the opposite orientation, then
have $$D(K^-)= d(\partial_+ K^-,\partial_-K^-)= \wt{lf}(\partial_+
K^-,\partial_- K^-)- lk(\partial_+ K^-,\partial_-K^-)=$$
$$\wt{lf}(\partial_+ K,\partial_-K)-
lk(\partial_+K,\partial_-K)=D(K),$$ by the bilinearity of
$\wt{lf}$ and of $lk.$
\end{proof}

The next theorem, which combines Theorems \ref{th2} and \ref{th3},
is the main result of the paper.

\begin{theorem}\label{main}
If $M$ satisfies condition (I) or (II) then there is an
isomorphism of $R$-algebras $$\Psi: {\cal S}_{\omega}(\pi_1(M),R,A)\to
{\cal S}(M,R,A),$$ such that $\Psi(x_g)=K\cdot A^{D(K)},$ where $K$
is an arbitrary framed unoriented knot which represents the
conjugacy class of $g^{\pm 1}\in \pi_1(M).$
\end{theorem}

We proceed the proof by a few introductory lemmas.

Let $\sim$ be the equivalence relation on $G$ such that
$x\sim y,\ x,y\in G,$ if and only if either $x$ or $x^{-1}$ is
conjugated to $y$ in $G.$ Each unoriented knot in 
$M$ determines a unique $\sim$-class in $\pi_1(M)$ which 
we denote by $<K>.$

\begin{lemma}\label{prop_of_d2}
\begin{enumerate}
\item[(1)] $A^{D(K)}K\in {\cal S}(M,R,A)$ depends on $<K>\in
\pi_1(M)/\sim$ only. In other words, $A^{D(K)}K=A^{D(K')}K',$ for
knots $K,K'$ representing the same $\sim$-class in $\pi_1(M).$
\item[(2)] Let $K_g, K_h$ be framed oriented knots representing the 
conjugacy classes of $g,h\in \pi_1(M).$
Assume that $K_h$ passes close to $K_g$ as presented on the
picture:\\ \centerline{\diag{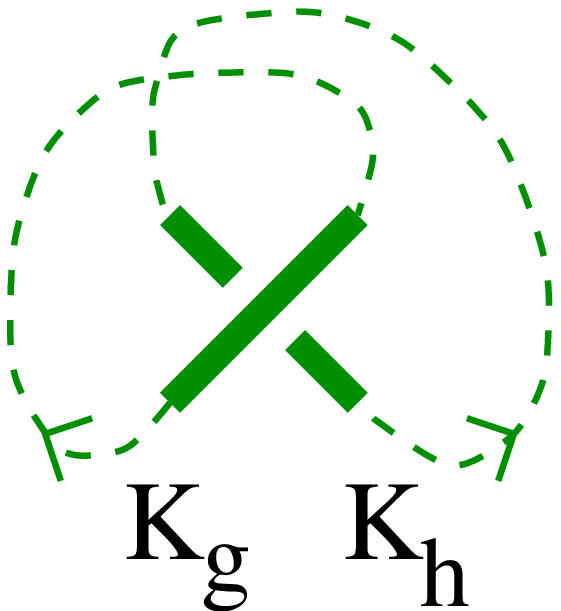}{2cm}} 
If $K_{gh}, K_{gh^{-1}}$ are knots obtained
by replacing \diag{3}{.6cm} in $K_g\cup K_h$ by \diag{4}{.6cm} and
\diag{5}{.6cm} respectively\footnote{The notation 
$K_{gh}, K_{gh^{-1}}$ is justified
by the fact that $K_{gh}$ and $K_{hg^{-1}}$ with the orientation
showed above represent $gh\in \pi_1(M)$ and $gh^{-1}\in \pi_1(M)$
respectively.},\\ \centerline{\diag{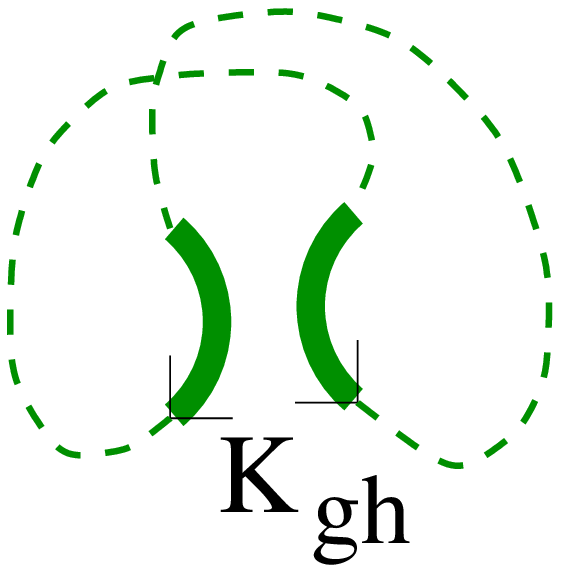}{2cm},\quad
\diag{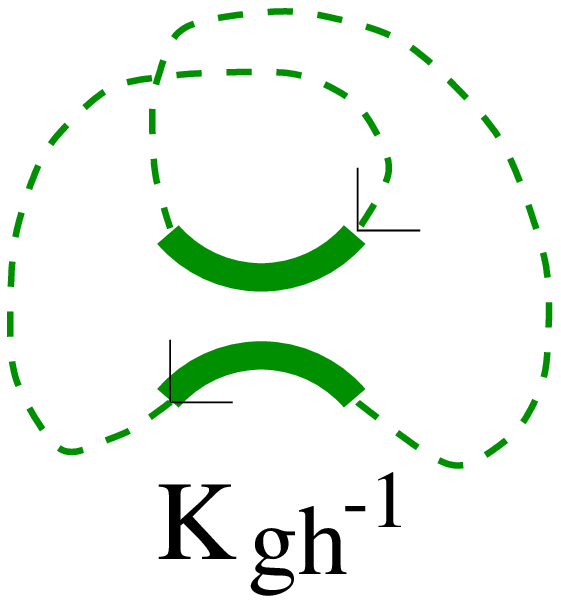}{2cm}} 
then
\begin{equation}
D(K_g)+D(K_h)+2d(K_g,K_h)=D(K_{gh})-1+\omega(g,h)
\label{D-identity-1}
\end{equation}
\begin{equation}
D(K_g)+D(K_h)+2d(K_g,K_h)=D(K_{gh^{-1}})+1-\omega(g,h).
\label{D-identity-2}
\end{equation}
\item[(3)] Let $K_g$ and $K_h$ be knots representing the conjugacy 
classes of $g,h\in \pi_1(M),$\\
\diag{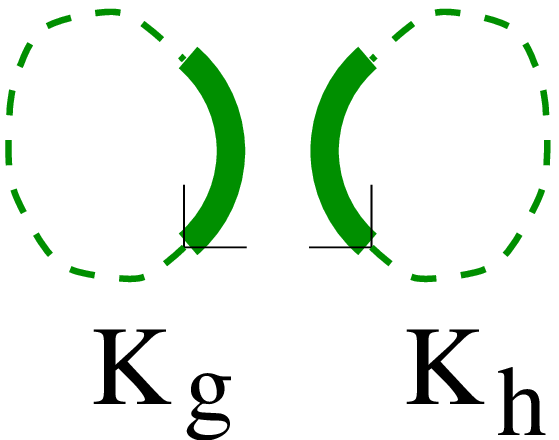}{1.5cm}.
\parbox{10cm}{If $K_{gh}, K_{gh^{-1}}$ represent knots obtained by replacing
\diag{4}{.6cm} by \diag{3}{.6cm} and
\diag{5}{.6cm} respectively, then}
\begin{equation}
D(K_g)+D(K_h)+2d(K_g,K_h)-D(K_{gh})-1=\omega(g,h)
\label{D-identity-3}
\end{equation}
and \begin{equation}
D(K_g)+D(K_h)+2d(K_g,K_h)-D(K_{gh^{-1}})-1=\omega(g,h).
\label{D-identity-4}
\end{equation}
\end{enumerate}
\end{lemma}

\begin{proof}
\begin{enumerate}
\item Two framed oriented knots represent the same $\sim$-class in
$\pi_1(M),$ if one can be obtained from the other by isotopy, crossing
changes, or by the reversal of orientation. By Lemma
\ref{d_no_orient}, the reversal of orientation of a knot $K$ does
not change $A^{D(K)}K\in {\cal S}(M,R,A).$ Suppose that framed
knots $K,$ $K'$ differ by a crossing change. Then the
$2$-component links, $\partial K$ and $\partial K',$ differ by two
crossing changes of the same sign. Hence, $$D(K')= d(\partial_+
K',\partial_-K')= \wt{lf}(\partial_+ K',\partial_- K')-
lk(\partial_+ K',\partial_-K')=$$ $$\wt{lf}(\partial_+
K,\partial_-K)-lk(\partial_+K,\partial_-K)\pm 2= D(K)\pm 2.$$
Hence, by (\ref{e4}), $$A^{D(K')}K'=A^{D(K)\pm 2}(-K)=
A^{D(K)}K.$$
\item By definition of $D,$ the left side of (\ref{D-identity-1}) is
$$\wt{lf}([K_g],[K_g])+\wt{lf}([K_h],[K_h])+
2\wt{lf}([K_g],[K_h])-lk(\partial_+ K_g,\partial_- K_g)+$$
$$-lk(\partial_+ K_h,\partial_- K_h)-2lk(K_g,K_h)=$$
$$\wt{lf}([K_g]+[K_h],[K_g]+[K_h])+
(\wt{lf}([K_g],[K_h])-\wt{lf}([K_h],[K_g]))+$$
$$-lk(\partial_+ K_g\cup \partial_+ K_h,\partial_- K_g\cup \partial_-
K_h).$$ 
By Lemma \ref{prop_of_d}(7),
$\wt{lf}([K_g],[K_h])-\wt{lf}([K_h],[K_g])=\omega(g,h).$
Since $[K_g+K_h]=[K_{gh}]$ in $H_1(M,\Z),$ the above expression
is equal to 
$$\wt{lf}([K_{gh}],[K_{gh}])+\omega(g,h)-
lk(\partial_+ K_g\cup \partial_+ K_h,\partial_- K_g\cup \partial_-
K_h).$$
Hence
\begin{equation}
\begin{array}{c}
D(K_g)+D(K_h)+2d(K_g,K_h)=\\
\wt{lf}([K_{gh}],[K_{gh}])+\omega(g,h)-
lk(\partial_+ K_g\cup \partial_+ K_h,\partial_- K_g\cup \partial_-
K_h).\label{long-e-1}
\end{array}
\end{equation}
$lk(\partial_+ K_g\cup \partial_+ K_h,\partial_- K_g\cup \partial_-
K_h)$ is the linking number between the components labeled by $1$ and
$2$ in \diag{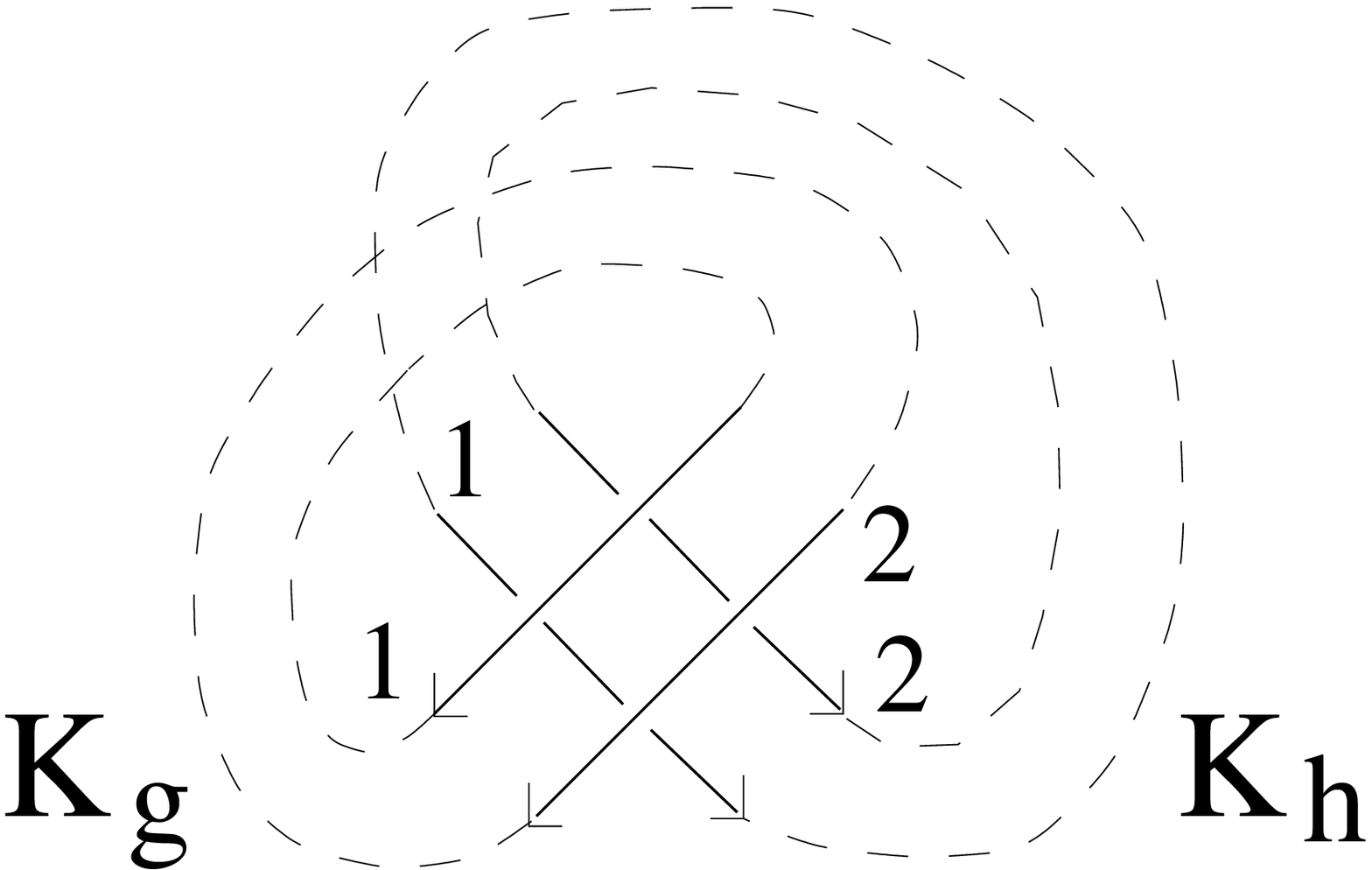}{1.8cm}.
Hence, it is equal to $1+lk(L_1,L_2),$ where $L_i,$ for $i=1,2,$ denotes the
components labeled by $i$ in \diag{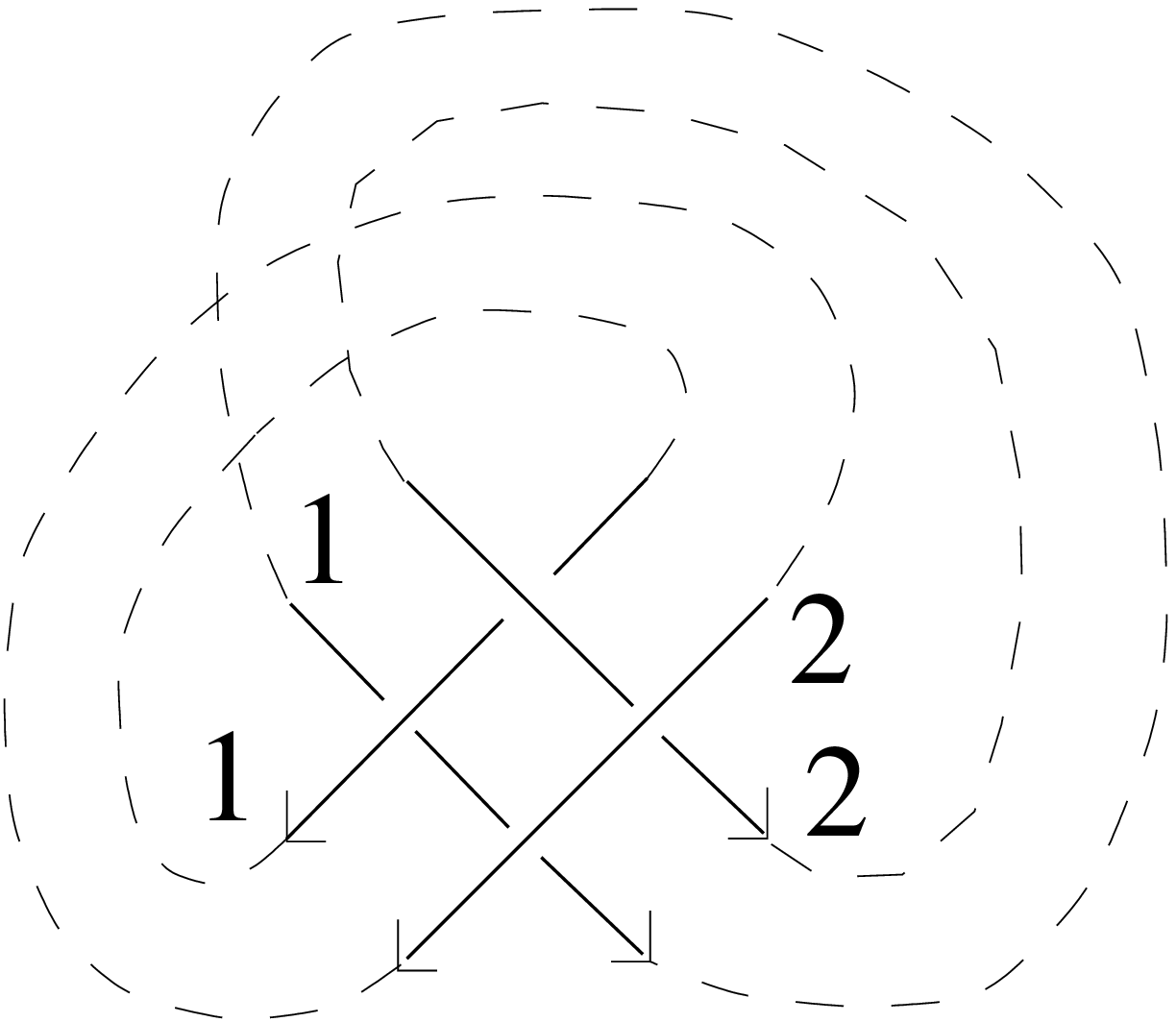}{1.8cm}.
Smoothing crossings in $L_1$ and in $L_2$
does not change $lk(L_1,L_2).$ Therefore, the above expression is equal to
$1+lk(L_1',L_2'),$ where $L_i'$ denotes the sublink
composed of components labeled by $i$ in \diag{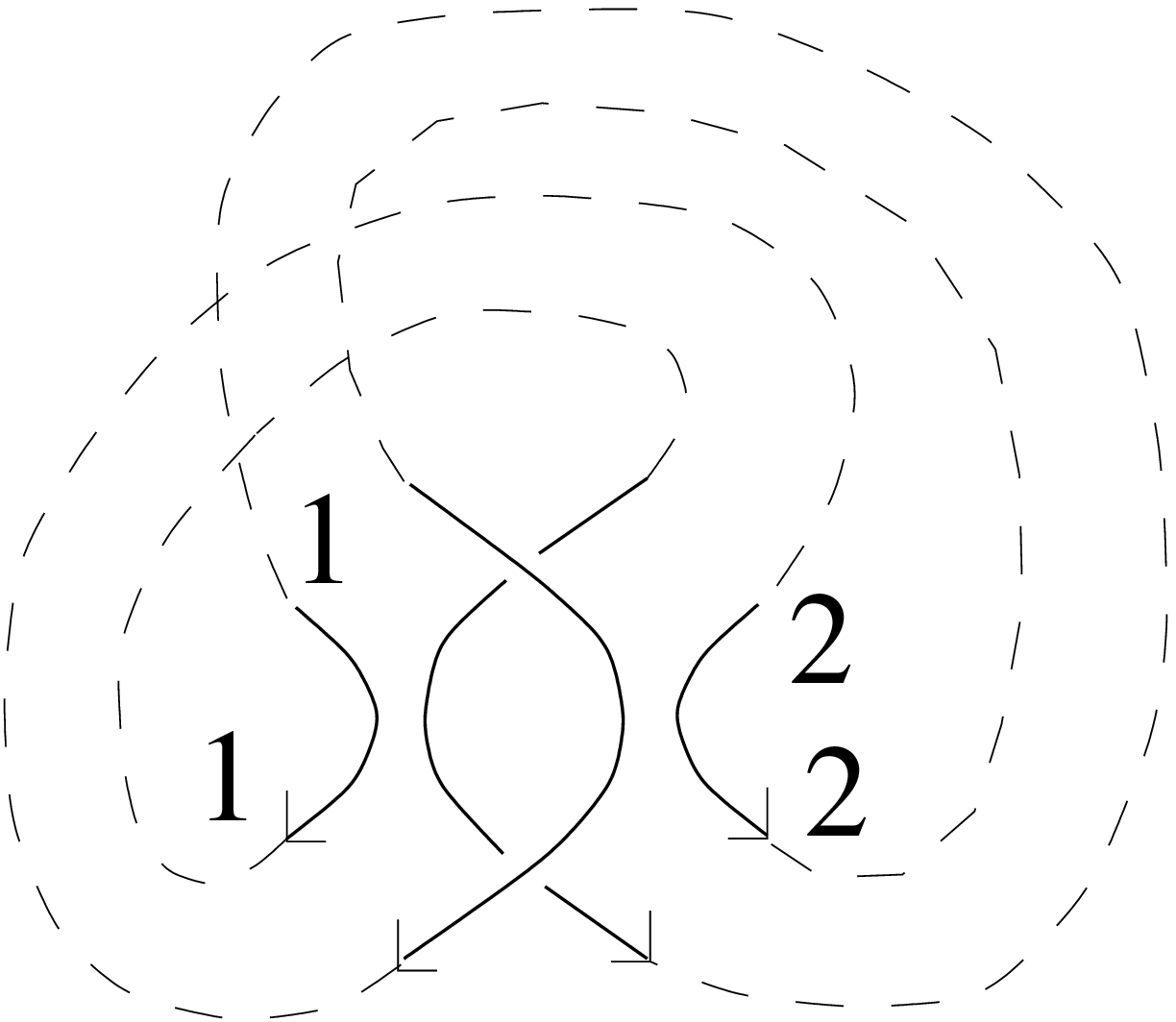}{1.8cm}.
Hence, we have proved that
$$lk(\partial_+ K_g\cup \partial_+ K_h,\partial_- K_g\cup \partial_-
K_h)=lk(\partial_+ K_{gh},\partial_- K_{gh})+1.$$
By applying this equation to (\ref{long-e-1}) we get (\ref{D-identity-1}):
$$D(K_g)+D(K_h)+2d(K_g,K_h)=$$
$$\wt{lf}([K_{gh}],[K_{gh}])+\omega(g,h)-
lk(\partial_+ K_{gh},\partial_- K_{gh})-1=$$
$$D(K_{gh})-1+\omega(g,h).$$
We prove identities (\ref{D-identity-2}), (\ref{D-identity-3}), and
(\ref{D-identity-4}) in the same way.
\end{enumerate}
\end{proof}

\begin{lemma}
If $M$ satisfies condition (I) or (II) then there is a
unique homomorphism of $R$-algebras $$\Psi: {\cal S}_{\omega}(\pi_1(M),R,A)\to
{\cal S}(M,R,A),$$ such that $\Psi(x_g)=K\cdot A^{D(K)},$ where $K$
is an arbitrary framed unoriented knot which represents the
conjugacy class of $g^{\pm 1}\in \pi_1(M).$
\end{lemma}

\begin{proof}
Let $$\Psi: R<x_g,\ g\in \pi_1(M)>\to {\cal S}(M,R, A)$$
be an an algebra homomorphism given by $\Psi(x_g)=A^{D(K)}K,$ 
where $K$ is an arbitrary framed unoriented knot which represents 
the conjugacy class of $g^{\pm 1}\in \pi_1(M).$
By Lemma \ref{prop_of_d2}(1), $\Psi(x_g)$ is well defied for each
$g\in G$ and $\Psi(x_{hgh^{-1}})=\Psi(x_g).$
We also have $\Psi(x_e-2)=0,$ since 
$\Psi(x_e)=A^{D(\diag{0}{.4cm})}\diag{0}{.5cm}=-(A^2+A^{-2})=2.$
Therefore, in order to show that $\Psi$ can be factored to a
homomorphism ${\cal S}_{\omega}(\pi_1(M),R,A)\to {\cal S}(M,R,A)$ we need
to prove that
\begin{equation}
\Psi(x_g)\Psi(x_h)=A^{\omega(g,h)}\Psi(x_{gh})+A^{-\omega(g,h)}
\Psi(x_{gh^{-1}}).\label{e9}
\end{equation}
Let $K_g$ and $K_h$ be two disjoint oriented framed knots
representing $g,h\in \pi_1(M).$ We may assume (by isotoping $K_h$ if
necessary) that $K_h$ passes very close to $K_g.$ Let $K_{gh},
K_{gh^{-1}}$ represent knots, which as in Lemma
\ref{prop_of_d2}(2), are obtained by replacing \diag{3}{.6cm} in
$K_g\cup K_h$ by \diag{4}{.6cm} and \diag{5}{.6cm} respectively.
Now equation (\ref{e9}) takes the form
$$\diag{8}{1.5cm}A^{D(K_g)}A^{D(K_h)}
(-1)^{d(K_g,K_h)}=A^{\omega(g,h)}A^{D(K_{gh})}K_{gh}+
A^{-\omega(g,h)}A^{D(K_{gh^{-1}})}K_{gh^{-1}}.$$ By Lemma
\ref{prop_of_d2}(2) this equation reduces to the skein equation:
$$\diag{3}{.6cm}=A\diag{4}{.6cm}+A^{-1}\diag{5}{.6cm}.$$
Therefore, we proved that $\Psi$ induces an $R$-algebra 
homomorphism $\Psi: {\cal S}_{\omega}(\pi_1(M),R,A)\to {\cal
S}(M,R,A).$
\end{proof}

The argument showing that $\Psi$ is an isomorphism will require
the following lemma.

\begin{lemma}\label{prop_of_S(G)} For any  $x_g, x_h\in {\cal
S}_{\omega}(G,R,A)$ we have:
\begin{enumerate}
\item $x_g=x_{g^{-1}}$
\item $x_{gh}=x_{hg}$
\item $x_hx_g=(-1)^{\omega(g,h)}x_gx_h$
\end{enumerate}
\end{lemma}

\begin{proof}
\begin{enumerate}
\item By substituting $g=e$ to (\ref{e3}), we get $2x_h=x_h+x_{h^{-1}}$ and,
therefore, $x_h=x_{h^{-1}}.$
\item By (\ref{e3.5}), $x_{gh}=x_{g(hg)g^{-1}}=x_{hg}.$
\item Observe that $x_hx_g= A^{\omega(h,g)}x_{hg}+A^{-\omega(h,g)}
x_{hg^{-1}},$ and $x_{gh}=x_{hg},$
$x_{hg^{-1}}=x_{(hg^{-1})^{-1}}= x_{gh^{-1}}.$ Since
$\omega(\cdot,\cdot)$ is skew-symmetric, we get
$$x_hx_g=A^{-\omega(g,h)}x_{gh}+A^{\omega(g,h)}x_{gh^{-1}}=$$
$$A^{2\omega(g,h)}\left(A^{\omega(g,h)}x_{gh}+A^{-\omega(g,h)}x_{gh^{-1}}
\right)=(-1)^{\omega(g,h)}x_gx_h.$$ We use here the fact that
$A^4=1$ and, therefore, $A^{3\omega(g,h)}=A^{-\omega(g,h)}.$
\end{enumerate}
\end{proof}

We are going to complete the proof of Theorem \ref{main} by
constructing an inverse homomorphism to $\Psi.$ This will imply
that $\Psi$ is an isomorphism.

Let $$\Phi: R{\cal L}(M)\to {\cal S}_{\omega}(\pi_1(M),R,A)$$ be a
homomorphism of $R$-modules which sends a link $L$ composed of
connected components $K_1,...,K_n$ to
$x_{g_1}x_{g_2}...x_{g_n}A^{-\sum_i D(K_i)-2\sum_{i<j}
d(K_i,K_j)},$ where $g_1,g_2,...,g_n$ are the elements of $\pi_1(M)$
(represented up to $\sim$-equivalence) by $K_1,K_2,...,K_n.$
Observe, that by (\ref{e3.5}) and Lemma \ref{prop_of_S(G)},
$x_g=x_h$ for any $g\sim h,$ $g,h\in \pi_1(M).$ Therefore the elements
$x_{g_1},...,x_{g_n}$ are uniquely defined. Furthermore, $D(K_i)$
and $d(K_i,K_j)$ mod $2$ do not depend on the orientations of $K_i$ and
$K_j.$ Therefore the following lemma shows that $\Phi$ is well
defined.

\begin{lemma}
$x_{g_1}x_{g_2}...x_{g_n}A^{-\sum_i D(K_i)-2\sum_{i<j}
d(K_i,K_j)}\in {\cal S}_\omega(\pi_1(M),R,A)$ does not depend on the manner
in which the connected components of $L$ are enumerated.
\end{lemma}

\begin{proof}
We need to show that if $\sigma$ is a permutation on $n$ symbols
then
$$x_{g_{\sigma(1)}}x_{g_{\sigma(2)}}...x_{g_{\sigma(n)}}A^{-\sum_i
D(K_{\sigma(i)})-2\sum_{i<j} d(K_{\sigma(i)},K_{\sigma(j)})}=$$
$$x_{g_1}x_{g_2}...x_{g_n}A^{-\sum_i D(K_i)-2\sum_{i<j}
d(K_i,K_j)}.$$ Since each permutation is a product of
transpositions, it is enough to assume that $\sigma$ is itself a
transposition, $\sigma=(s,s+1).$ Notice that $\sum_i D(K_i)=\sum_i
D(K_{\sigma(i)}).$ Observe also, that the sums $$\sum_{i<j}
d(K_{\sigma(i)},K_{\sigma(j)}) {\rm \ and\ } \sum_{i<j}
d(K_i,K_j)$$ differ only by one term, that is, the first sum
contains $d(K_{s+1},K_s),$ and the second sum contains
$d(K_s,K_{s+1}).$ Therefore the above equality reduces to
$$x_{g_{s+1}}x_{g_s}A^{-2d(K_{s+1},K_s)}=x_{g_s}x_{g_{s+1}}
A^{-2d(K_s,K_{s+1})},$$ which follows from Lemmas
\ref{prop_of_d}(7) and \ref{prop_of_S(G)}(3).
\end{proof}

\begin{lemma} The submodule
${\cal S}\subset R{\cal L}(M)$ generated by expressions
(\ref{e1}) and (\ref{e2}) is contained in the kernel of $\Phi.$
\end{lemma}

\begin{proof}
We need to show that expressions (\ref{e1}) and (\ref{e2})
belong to the kernel of $\Phi.$ 
Let us look at (\ref{e2}) first:
if a link $L=K_1\cup ...\cup K_n$ is unlinked with \diag{0}{.5cm}
then $\Phi(L\cup \diag{0}{.5cm})=\Phi(L)\cdot x_e$ since
$D(\diag{0}{.5cm})=0$ and $d(\diag{0}{.5cm}, L)=0$ by Lemma 
\ref{prop_of_d}(6). Hence, skein expression (\ref{e2}) lies in $Ker\,
\Phi.$
Now, we will show that (\ref{e1}) also belongs to $Ker\, \Phi$
by analyzing two cases, depending on how the ends of the tangle
\diag{3}{.6cm} are connected in the ambient space:\\
(1) Suppose that \diag{3}{.6cm} in (\ref{e1}) represents a link of 
$n$ components,
$K_{g_1}\cup ... \cup K_{g_n},$ whose the first two components
are depicted in the diagram below:\\
\centerline{\diag{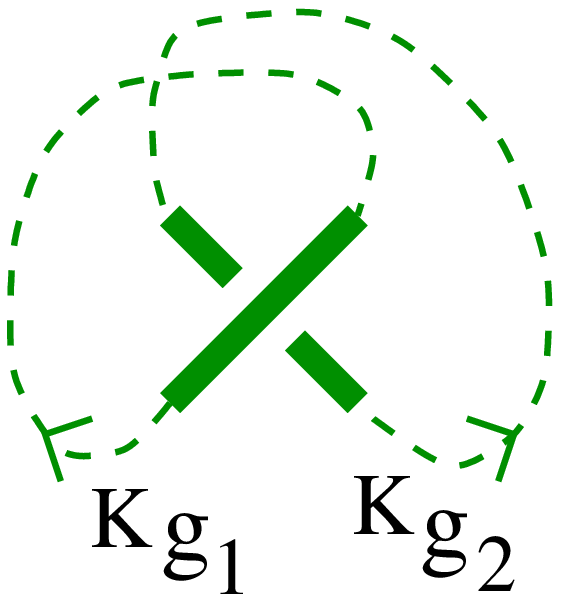}{2cm}.}
Let $K_{g_1g_2}$ and $K_{g_1g_2^{-1}}$ be as in Lemma
\ref{prop_of_d2}(2),\\ \centerline{\diag{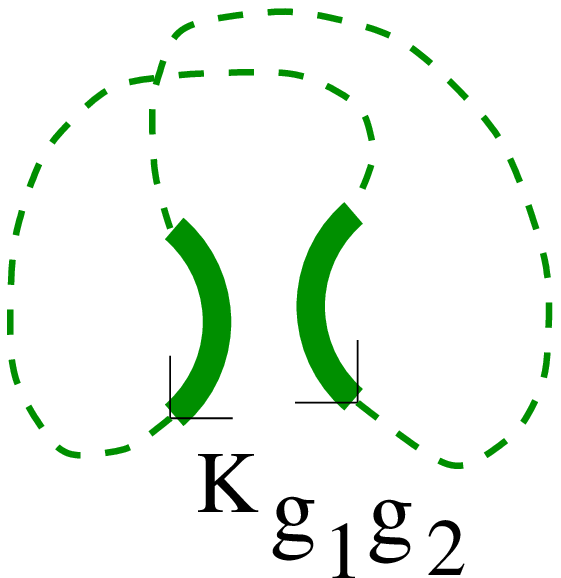}{2cm},\quad
\diag{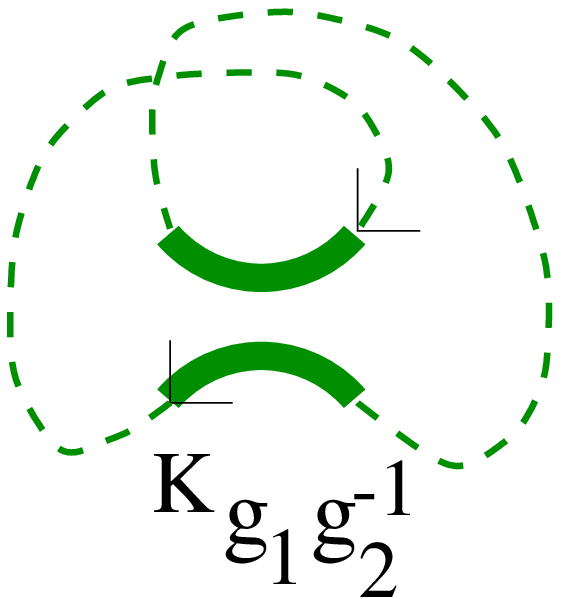}{2cm}.}
Then we need to show that
\begin{equation}
\begin{array}{c}
\Phi(K_{g_1}\cup K_{g_2}\cup ... \cup K_{g_n})=\\
A\Phi(K_{g_1g_2}\cup K_{g_3}\cup ... \cup K_{g_n})+
A^{-1}\Phi(K_{g_1g_2^{-1}}\cup K_{g_3}\cup ... \cup K_{g_n}).
\end{array}
\label{skein-phi-1}
\end{equation}
The above equation, can be rewritten as
$$x_{g_1}x_{g_2}...x_{g_n}A^{-\sum_i D(K_{g_i})-2\sum_{i<j}
d(K_{g_i},K_{g_j})}=$$
$$Ax_{g_1g_2}x_{g_3}...x_{g_n}A^{-D(K_{g_1g_2})-
\sum_{i>2} D(K_{g_i})-2\sum_{i>2} d(K_{g_1g_2},K_{g_i})- 2\sum_{2<i<j}
d(K_{g_i},K_{g_j})}+$$
$$A^{-1}x_{g_1g_2^{-1}}x_{g_3}...x_{g_n}
A^{-D(K_{g_1g_2^{-1}})- \sum_{i>2} D(K_{g_i})-2\sum_{i>2} 
d(K_{g_1g_2^{-1}},K_{g_i})- 2\sum_{2<i<j}
d(K_{g_i},K_{g_j})}.$$
By Lemma \ref{prop_of_d}(5), $$d(K_{g_1},K_{g_i})+d(K_{g_2},K_{g_i})=
d(K_{g_1g_2},K_{g_i})=d(K_{g_1g_2^{-1}},K_{g_i}),$$ for $i>2.$
Therefore the above equation simplifies to
$$x_{g_1}x_{g_2}A^{-D(K_{g_1})-D(K_{g_2})-2d(K_{g_1},K_{g_2})}=
Ax_{g_1g_2}A^{-D(K_{g_1g_2})}
+A^{-1}x_{g_1g_2^{-1}}A^{-D(K_{g_1g_2^{-1}})}.$$
By (\ref{D-identity-1}) and (\ref{D-identity-2}), the above equation
reduces to (\ref{e3}). This completes the proof of
(\ref{skein-phi-1}).\\
(2) The second possibility is that \diag{3}{.6cm} in (\ref{e1}) 
represents a link of $n$ components, 
$K_{g_1g_2}\cup K_{g_3}... \cup K_{g_n},$
whose the first component
is depicted in the diagram below:\\
\centerline{\diag{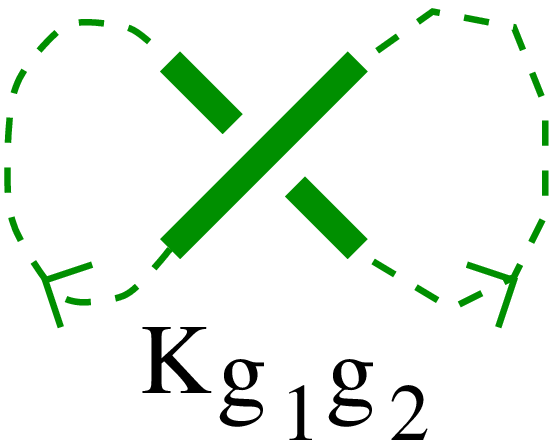}{1.5cm}.}
Denote \diag{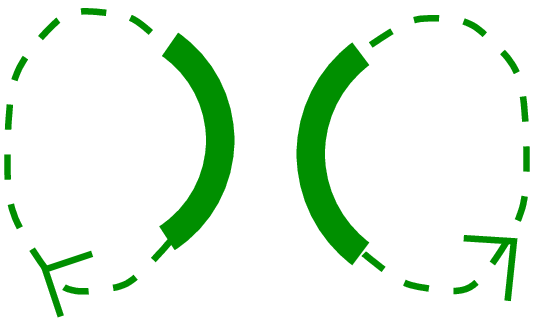}{.8cm} and \diag{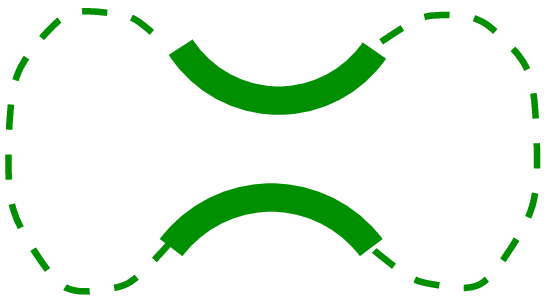}{.8cm} by
$K_{g_1}\cup K_{g_2}$ and $K_{g_1g_2^{-1}}$ respectively.
By the method used for proving (1), we prove that
$$\Phi(K_{g_1g_2}\cup K_{g_3}\cup ... \cup K_{g_n})=$$
$$A\Phi(K_{g_1}K_{g_2}\cup ... \cup K_{g_n})+
A^{-1}\Phi(K_{g_1g_2^{-1}}\cup K_{g_3}\cup ... \cup K_{g_n}).$$
We use identities (\ref{D-identity-3}) and 
(\ref{D-identity-4}) in the proof.
\end{proof}

\begin{lemma}
$\Phi$ is a homomorphism of algebras.
\end{lemma}

\begin{proof}
Since ${\cal S}(M,R,A)$ is spanned by knots, it is enough to show that
$\Phi(K_1\cdot K_2)$ is equal to $\Phi(K_1)\cdot \Phi(K_2)$ for any
knots $K_1,K_2.$
We have
$$\Phi(K_1\cdot K_2)=\Phi(K_1\cup K_2\cdot
(-1)^{d(K_1,K_2)})=$$
$$x_{g_1}x_{g_2}A^{-D(K_1)-D(K_2)-2d(K_1,K_2)}(-1)^{d(K_1,K_2)}=$$
$$\Phi(K_1)\cdot \Phi(K_2).$$
\end{proof}

\noindent{\bf Proof of Theorem \ref{main}:}
We have $$\Phi(\Psi(x_g))=\Phi(K\cdot
A^{D(K)})=x_gA^{-D(K)}A^{D(K)}=x_g$$
and, similarly, $\Psi(\Phi(K))=K,$ for any $g\in \pi_1(M)$ and any
knot $K$ in $M.$ Therefore, $\Phi$ is the inverse of $\Psi,$ and
hence $\Psi$ is an isomorphism.
\Box

%
\end{document}